% Logic Eprints
%Submitted 1740 Sun Jul 17, 1994 by: jobr@michelangelo.mathematik.uni-tuebingen.de (joerg brendle)
%logic/brendle/mupe.tex
%

\magnification=\magstep1

\def\Bool{[\![}
\def\Boor{]\!]}

%Das Folgende sollte die "Blackboard bold font" \Bbb definieren, und 
%auch  die "Fraktur" font, \frak:
% Beispiel:    $ {\frak c} = |{\Bbb R}| $

\def\hexnumber#1{\ifcase#1 0\or1\or2\or3\or4\or5\or6\or7\or8\or9\or
	A\or B\or C\or D\or E\or F\fi }

%  The following lines establish the use of the Euler Fraktur font.
\font\teneuf=eufm10
\font\seveneuf=eufm7
\font\fiveeuf=eufm5
\newfam\euffam
\textfont\euffam=\teneuf
\scriptfont\euffam=\seveneuf
\scriptscriptfont\euffam=\fiveeuf

%  End definition of Euler Fraktur font.

\font\tenmsx=msam10
\font\sevenmsx=msam7
\font\fivemsx=msam5
\font\tenmsy=msbm10
\font\sevenmsy=msbm7
\font\fivemsy=msbm5
\newfam\msxfam
\newfam\msyfam
\textfont\msxfam=\tenmsx  \scriptfont\msxfam=\sevenmsx
  \scriptscriptfont\msxfam=\fivemsx
\textfont\msyfam=\tenmsy  \scriptfont\msyfam=\sevenmsy
  \scriptscriptfont\msyfam=\fivemsy
\edef\msx{\hexnumber\msxfam}

\mathchardef\upharpoonright="0\msx16
\let\restriction=\upharpoonright
\def\Bbb#1{\tenmsy\fam\msyfam#1}

\def\restrict{{\restriction}}
\def\re{{\restriction}}

\def\Smallskip{\vskip1.4truecm}
\def\Bigskip{\vskip2.2truecm}

\def\qed{{\vcenter{\hrule height.4pt \hbox{\vrule width.4pt height5pt
 \kern5pt \vrule width.4pt} \hrule height.4pt}}}
\def\ok{\vbox{\hrule height 8pt width 8pt depth -7.4pt
    \hbox{\vrule width 0.6pt height 7.4pt \kern 7.4pt \vrule width 0.6pt height 7.4pt}
    \hrule height 0.6pt width 8pt}}
\def\nt{{\leq}\kern-1.5pt \vrule height 6.5pt width.8pt depth-0.5pt \kern 1pt}
\def\sd{{\times}\kern-2pt \vrule height 5pt width.6pt depth0pt \kern1pt}
\def\zp#1{{\hochss Y}\kern-3pt$_{#1}$\kern-1pt}

\def\AA{{\Bbb A}}

\def\BB{{\Bbb B}}
\def\CC{{\Bbb C}}

\def\PP{{\Bbb P}}
\def\QQ{{\Bbb Q}}
\def\RR{{\Bbb R}}
\def\SS{{\Bbb S}}

\def\A{{\cal A}}
\def\B{{\cal B}}

\def\L{{\cal L}}

\def\sm{{\smallskip}}
\def\ce#1{{\centerline{#1}}}
\def\no{{\noindent}}
\def\la{{\langle}}
\def\ra{{\rangle}}
\def\sub{\subseteq}

\def\ha{{\hat{\;}}}
\def\em{{\emptyset}}
\def\sem{\setminus}
\def\omom{{\omega^\omega}}
\def\omlom{{\omega^{<\omega}}}

\def\twoom{{2^\omega}}
\def\twolom{{2^{<\omega}}}

\def\Loriar{\Longrightarrow}
\def\Loleriar{\Longleftrightarrow}
\font\small=cmr8 scaled\magstep0
\font\smalli=cmti8 scaled\magstep0
\font\capit=cmcsc10 scaled\magstep0
\font\capitg=cmcsc10 scaled\magstep1
\font\dunh=cmdunh10 scaled\magstep0
\font\dunhg=cmdunh10 scaled\magstep1
\font\dunhgg=cmdunh10 scaled\magstep2

\font\sanse=cmss10 scaled\magstep0

\font\bolds=cmssdc10 scaled\magstep0

\overfullrule=0pt
\openup1.5\jot

\ce{}
\Smallskip
\ce{\dunhgg   Mutually Generics and Perfect Free Subsets }
\footnote{}{{\openup-6pt {\small {\smalli
1991 Mathematics subject classification.} 03E05 03E15 03E35 03E40
     \par
{\smalli Key words and phrases.} perfect set, perfect tree, 
superperfect tree, forcing, Cohen real, random
real, Sacks real, free subset, Borel function, measurable function.
\endgraf}}}
\Bigskip
\ce{\capitg J\"org Brendle\footnote{$^\star$}
{{\small Supported by DFG--grant Nr. Br 1420/1--1.}}}
\Smallskip
\no Mathematisches Institut der Universit\"at
T\"ubingen, Auf der Morgenstelle 10, 72076
T\"ubingen, Germany; email: {\sanse
jobr@michelangelo.mathematik.uni--tuebingen.de}
\Bigskip
\ce{\capit Abstract}
\bigskip
\no We present two ways of adjoining a perfect set of mutually
random reals to a model $V$ of $ZFC$. We also investigate the
existence of perfect free subsets for projective functions
$f : (\omom)^n \to \omom$.
\vfill\eject

{\dunhg Introduction}
\Smallskip

This work is about various aspects of perfect sets of reals
the elements of which are independent of each other in a
strong sense. As is usual in the area, when talking about
the reals, we shall mean the Cantor space $\twoom$ or
the Baire space $\omom$ rather than the real line $\RR$ itself.
\par

Recall that a tree $T \sub \twolom$ is said to be {\it
perfect} if for all $\tau\in T$ there is $\sigma \supseteq
\tau$ with $\sigma \ha\la 0\ra \in T$ and $\sigma\ha\la 1\ra \in T$.
For a perfect tree $T$ we let $[T] := \{ f \in\twoom ; \;
\forall n \in\omega \; (f\re n \in T ) \}$ denote the
set of its branches. Then $[T]$ is a perfect set  (in the
topology of $\twoom$). Conversely, given a perfect set
$S \sub\twoom$ there is a perfect tree $T \sub\twolom$ such that
$[T] = S$. This allows us to confuse perfect sets and
perfect trees in the sequel. ---
It is well--known that adding one Cohen real $c$ to a model $V$
of $ZFC$ in fact adds a perfect set of {\it mutually
generic Cohen reals}. This means that in the generic extension $V[c]$,
there is a perfect tree $T \sub \twolom$ such that given $n
\in\omega$ and distinct $x_0 , ... , x_n \in [T]$, $x_n$ is
Cohen--generic over $V [x_0 , ... , x_{n-1} ]$ (see 1.1
for details). In section 2, we shall describe situations
in which a perfect set of {\it mutually generic random
reals} is adjoined. This notion is defined in an exactly
similar fashion, with {\it Cohen} replaced by {\it random}.
Improving an earlier result of Bartoszy\'nski and Judah
[BJ 1, Theorem 2.7], we show:
\sm

{\bolds Theorem A.} {\it Let $V \sub W$ be models of
$ZFC$. Assume there is a dominating real $d$ over $V$ in $W$,
and $r$ is random over $W$; then there is a perfect set
$T$ of mutually random reals over $V$ in $W [r]$. }
\sm

\no Here, a real $d\in\omom\cap W$ is called a {\it dominating
real over} $V$, for models $V\sub W$ of $ZFC$, if for all
$f\in \omom\cap V$, $f(n) < d(n)$ holds for all but finitely
many $n\in\omega$. ---
These results shed new light on old theorems of Mycielski's
(see [My 1], [My 2]) who proved that given a sequence of Borel
null (meager, respectively) sets $\la B_n \sub \RR^{k(n)} ; \;
n\in\omega\ra$, where $k(n) \in\omega$, there is a perfect set
$A \sub \RR$ such that $A^{k(n)} \cap B_n \sub \{ (x_0 , ... ,
x_{k(n) - 1} ) ; \; \exists i \neq j \; (x_i = x_j ) \}$ for all
$n\in\omega$. Clearly, we can get such a set $A$ by taking
a generic for the p.o. of Theorem A (of 1.1, respectively)
over a countable model of $ZFC$ containing the Borel codes of the 
$B_n$'s.
\par

To appreciate the other problem we are dealing with,
assume we are given a function $f :
\RR^n \to \RR$. A set $A \sub \RR$ is {\it free for} $f$ if
for all $\la a_0 , ... , a_{n-1} \ra \in A^n$, we have
$f (\la a_0 , ... , a_{n-1} \ra ) \in \{ a_0 , ... , a_{n-1} \}
\cup (\RR \sem A)$. A function $f : \RR^m \to \RR$ is
{\it $\Delta^1_n$--measurable} (a {\it $\Delta^1_n$--function}
for short) if for all open $U \sub \RR$, $f^{-1} (U)$ is
a $\Delta^1_n$--set. The {\it Borel functions} are just the 
$\Delta_1^1$--functions. ---
Mildenberger [Mi, Theorem 1.1] proved that
if every $\Delta^1_n$--set has
the property of Baire, then every $\Delta^1_n$--function
$f : (\twoom)^m \to \twoom$ has a perfect free subset.
This can also be shown using Mycielski's Theorem [My 1].
Similarly, with the aid of [My 2], one gets that if
every $\Delta^1_n$--set is measurable, then every $\Delta^1_n$--function
$f : (\twoom)^m \to \twoom$ has a perfect free subset
(see 1.4 for details).
\par
We shall henceforth write $\Delta^1_n ({\cal B})$
($\Sigma^1_n ({\cal B}), \Delta^1_n (\L),
\Sigma^1_n (\L)$) for {\it every $\Delta^1_n$($\Sigma^1_n$)--set
has the Baire property} ({\it is Lebesgue measurable}), 
and $\Delta^1_n ({\cal PFS}_m)$ for
{\it every $\Delta^1_n$--function $f : (\twoom)^m \to \twoom$ has a perfect
free subset}. $\Delta^1_n ({\cal PFS})$ abbreviates $\forall m
\; (\Delta^1_n ({\cal PFS}_m))$. The above discussion 
shows that $\Delta^1_n ({\cal PFS})$ is weaker than both
$\Delta^1_n (\L)$ and $\Delta^1_n (\B)$. In $\S$ 3 we prove that
it is consistently strictly weaker for $n \geq 3$.
\sm

{\bolds Theorem B.} {\it $\forall n \; (\Delta^1_n ({\cal PFS})) +
\neg \Sigma^1_2 ({\cal B}) + \neg \Delta^1_2 (\L)$ is consistent with $ZFC$.}
\sm

\no We do not know whether $\forall n \; (\Delta_n^1 ({\cal PFS}))
+ \neg \Delta^1_2 ({\cal B})$  is consistent, but
conjecture the answer to be positive.
\par
On the other hand, $\Delta^1_2 ({\cal PFS})$ is not a Theorem of
$ZFC$, for Mildenberger [Mi, Theorem 1.3] also proved that $V=L$ implies that
$\Delta^1_2 ({\cal PFS}_1)$ fails. Still, $\Delta^1_2 ({\cal PFS}_1)$ has a very weak
characterization. It is equivalent to $V \neq L[x]$ for all
reals $x$ (see 4.2). Using this we shall show (still in
$\S$ 4):
\sm
 
{\bolds Theorem C.} {\it $\Delta^1_2 ({\cal PFS}_1) + \neg
\Delta^1_2 ({\cal PFS}_2)$ is consistent with $ZFC$.}
\sm

\no This coincides with our intuition that 
it is much easier to construct a
``large" free subset for functions $f: \RR^m \to \RR$ in case
$m=1$ than in case $m \geq 2$. 
\par

A tree $T\sub\omlom$ is called {\it superperfect} if for all
$\sigma \in T$ there is $\tau \supseteq \sigma$ such that $\tau\ha
\la n\ra \in T$ for infinitely many $n\in\omega$. A set of
branches through a superperfect tree shall be called a {\it
superperfect set}. We conclude our considerations with
some results on superperfect trees in section 5.
In particular we will see that most of the results
on perfects sets cannot be extended to superperfect
sets.
\bigskip

{\bolds Notational Remarks.} Most of our notation should
follow set--theoretic convention (see e.g. [Je 1]). We just explain a few
things which may be less standard.
\par
Whenever we write $\Delta^1_n$, ... 
we mean the boldface version. $<_c$ denotes complete embedding
of Boolean algebras. A p.o. $\PP$ is {\it $\sigma$--linked}
if there are $P_n \sub \PP$ ($n\in\omega$) with $\bigcup_n P_n = \PP$ and 
for all $n\in\omega$ any two elements of $P_n$ are compatible.
For random forcing $\BB$, Cohen forcing $\CC$ and Sacks forcing $\SS$
we refer the reader to [Je 2, chapter I, section 3].
\par
$\forall^\infty n$ denotes {\it for almost all $n$}, while
$\exists^\infty n$ stands for {\it there are infinitely many $n$}.
Given $\sigma \in\twolom$, $[\sigma] := \{ x\in
\twoom ; \; \sigma \sub x\}$ is the clopen set determined by $\sigma$.
Given a superperfect tree $T \sub \omlom$, $split (T) =
\{ \sigma \in T ; \; \exists^\infty n \; (\sigma \ha \la n\ra
\in T ) \}$ is the set of {\it $\omega$--splitting nodes}
of $T$. Given models $V\sub W$ of $ZFC$, $f \in \omom\cap W$
is an {\it unbounded real over} $V$ if for all $g\in\omom\cap V
\; \exists^\infty n \; (g(n) < f(n))$; and $T\in W$ is a {\it perfect
tree of Cohen} ({\it random}, respectively) {\it reals over} $V$
if every branch of $[T]\cap W$ is Cohen (random, resp.) over $V$.
(This is a weakening of the notion a perfect tree of mutually Cohen (random, resp.)
reals.)
\Bigskip

{\dunhg $\S$ 1. Preliminary facts}
\Smallskip

We list here a few results which are fundamental for our
work. The first two may belong to set--theoretic folklore. However,
as we could not find appropriate references, we include
proofs.
\sm
{\capit 1.1. Lemma.} (Folklore) {\it Cohen forcing $\CC$ adds a perfect
set of mutually generic Cohen reals.}
\sm
{\it Proof.} Let $\PP = \{  (t,n) ; \; t \sub \twolom$ is a finite subtree
$\land\; n\in\omega\;\land\; t$ has height $n\;\land\;$ all branches
of $t$ have length $n \}$, ordered by end--extension; i.e.
$(t,n) \leq (s,m)$ if $t \supseteq s$ and $n \geq m$ and every
branch of $t$ extends a branch of $s$. $\PP$ is a countable
p.o. and thus forcing equivalent to $\CC$. It generically adds
a perfect tree $T$. We have to check that the branches of $T$
are mutually Cohen.
\par
To see this take $f_0 , ... , f_m\in [T]$ in the extension $V[G]$ via
$\PP$. Take a condition $(t,n) \in \PP \cap G$ in the ground model
$V$ such that there are distinct $\sigma_0 , ... , \sigma_m \in t$,
$| \sigma_i | = n - 1$ ($i \leq m$), with $\sigma_i \sub f_i$.
Let $A \sub [\twoom]^{m+1}$ be a nowhere dense set coded in $V$.
Note that given any $(s , \ell) \leq (t,n)$ in $\PP$, there is
$(r,k) \leq (s,\ell)$ such that: whenever $\sigma_0 ' , ... , \sigma_m '
\in r$, $| \sigma_i ' | = k-1$, satisfy $\sigma_i \sub \sigma_i '$, then 
$[\sigma_0 '] \times ... \times [\sigma_m '] \cap A = \em$
$(\star)$. By genericity, a condition $(r,k)$ with $(\star)$ must lie
in the filter $G$. Thus $\la f_0 , ... , f_m \ra \notin A$, and
$\la f_0 , ... , f_m \ra$ is generic for $\CC^{m+1}$; this entails, however,
that any $f_i$ is $\CC$--generic over the others. $\qed$
\sm
Note that to say that $T$ is perfect tree of mutually Cohen reals
over $V$ is absolute for $ZFC$--models $W \supseteq V$ with $T \in W$. 
The point is that this is equivalent to saying that for all $m$ the set
of sequences of length $m$ of distinct branches of $T$ has
empty intersection with every $m$--dimensional meager set in $V$.
The latter statement is absolute (see [Je 1, section 42]). A similar
remark applies to perfect sets of mutually random reals --- and
also to perfect sets of Cohen (random) reals.
\bigskip

Given a function $f : (\omom)^m \to \omom$, we denote
by $G(f)$ the {\it graph} of $f$, i.e.
$G(f) = \{ (x , f(x)) ; \; x \in (\omom)^m \}$.
\sm
{\capit 1.2. Lemma.} (Folklore) {\it The following are equivalent
for a function $f : (\omom)^m \to \omom$:
\par
\item{(i)} $f$ is a $\Delta^1_n$--function; \par
\item{(ii)} $f$ is $\Sigma^1_n$--measurable (i.e. inverse images
of open sets are $\Sigma^1_n$); \par
\item{(iii)} $f$ is $\Pi_n^1$--measurable; \par
\item{(iv)} $G(f)$ is $\Delta^1_n$; \par
\item{(v)} $G(f)$ is $\Sigma^1_n$. \par }
\sm
{\it Proof.} The equivalence of (i) thru (iii) and the
implication (iv) $\Loriar$ (v) are obvious. \par
(i) $\Loriar$ (iv). $(x,y) \in G(f) \Loleriar
\forall k \; (x \in f^{-1} ( [y \re k]))$. \par
(v) $\Loriar$ (ii). Let $\sigma \in\omlom$. Then
\sm
\ce{$x \in f^{-1} ([\sigma]) \Loleriar \exists y \;
(y \in [\sigma] \;\land\; (x,y) \in G(f))$. \hskip
2truecm $\qed$}
\sm
\sm
\no Note that to say that $G(f)$ is $\Pi_n^1$ is weaker
than (i) --- (v) in the Lemma (see [Mo, chapter 5] for the
construction of a $\Pi_1^1$--graph in $L$ with unpleasant
properties). Thus we get the following hierarchy
of real functions:
\sm
\ce{Borel functions $\sub \Pi_1^1$--graphs $\sub \Delta^1_2$--functions
$\sub \Pi^1_2$--graphs $\sub ...$}
\bigskip

We state again explicitly Mycielski's Theorems. For a proof
see either [My 1] and [My 2] or use Lemma 1.1 and Theorem A (as
remarked in the Introduction). We note that Mycielski's original
formulation of these results is somewhat different. However,
it is rather easy to prove his version from the present one,
and vice--versa.
\sm
{\bolds 1.3. Theorem.} (Mycielski) {\it (a) Let $\la B_n \sub (\twoom)^{k(n)}
; \; n\in\omega \ra$ be a sequence of meager sets, where
$k(n) \in \omega$; then there is a perfect tree $T \sub \twolom$
satisfying $[T]^{k(n)} \cap B_n \sub \{ (x_0 , ... , x_{k(n) - 1} ); 
\; \exists i \neq j \; (x_i = x_j) \}$ for all $n\in\omega$. \par
(b) Let $\la B_n \sub (\twoom)^{k(n)} ; \; n\in\omega\ra$ be
a sequence of null sets, where $k(n) \in\omega$; then there
is a perfect tree $T \sub \twolom$ satisfying $[T]^{k(n)}
\cap B_n \sub \{ (x_0 , ... , x_{k(n) - 1} ) ; \; \exists
i \neq j \; (x_i = x_j) \}$ for all $n\in\omega$. $\qed$}
\sm

{\bolds 1.4. Theorem.} (a) (Mildenberger [Mi]) {\it $\Delta^1_n ({\cal B})$
implies $\Delta^1_n ({\cal PFS})$ for all $n\in\omega$. Furthermore,
$\Sigma^1_n (\B)$ implies that every $\Pi^1_n$--graph has
a perfect free subset. \par}
(b) {\it $\Delta^1_n (\L)$ implies $\Delta^1_n ({\cal PFS})$ for
all $n\in\omega$. Furthermore, $\Sigma^1_n (\L)$ implies that every
$\Pi^1_n$--graph has a perfect free subset.}
\sm
{\it Proof.} All four cases are similar. Therefore we restrict
ourselves to proving $\Delta^1_n (\L) \Loriar \Delta^1_n ({\cal PFS})$.
\par
Let $f: (\omom)^m \to \omom$ be a $\Delta^1_n$--function.
By 1.2, $G(f)$ is $\Delta^1_n$; by assumption $G(f)$ is
measurable; and by Fubini's Theorem $G(f)$ must be null.
Next, given a partition $\A=\{ A_i ; \; i < k \}$ of $m$
for some $k < m$, we let $X_{\A} = \{ (x_0 , ... , x_{m-1}
) \in (\omom)^m ; \; \forall i < k \; \forall j , j' \in n \;
(j,j' \in A_i  \Loriar x_j = x_{j'} ) \}$. Let $G_{\A}$
be the graph of the restriction $f \re X_{\A}$. As before,
all $G_{\A}$ are null in the space $X_{\A} \times \omom$.
Apply Mycielski's Theorem (1.3 (b)) to $G(f)$ and all the
$G_{\A}$'s, and get a perfect tree $T \sub \omlom$.
$[T]$ is easily seen to be a perfect free subset. $\qed$
\sm

\no A consequence of this is that $\Pi^1_1$--graphs always have
perfect free subsets. Thus Mildenberger's Theorem [Mi, Theorem 1.3]
is best possible.
\Bigskip

{\dunhg $\S$ 2. Perfect sets of mutually generic random reals}
\Smallskip

Cicho\'n (see [BJ 1, 2.1]) observed that random forcing does
not add a perfect set of random reals. Shelah [BrJS, 3.1]
strengthened this result by showing that the existence
of a perfect set of random reals over $V$ in a model $W \supset V$
implies that either $W$ contains a real eventually dominating
the reals of $V$ or it contains a null set $N$ with $\twoom \cap V
\sub N$. On the other hand, Bartoszy\'nski and Judah [BJ 1,
2.7] proved that if $W \supset V$ contains a real $d$
dominating the reals of $V$ and $r$ is random over $W$, then
$W[r]$ contains a perfect set of random reals over $V$.
Theorem A strengthens this; its proof follows
along the same lines. (However, we think that our
argument is somewhat shorter and more straightforward.)
\bigskip

{\it 2.1. Proof of Theorem A.} 
Using the dominating function $d$, define
the following sequence (in $W$):
\sm
\ce{$d_0 = 0 , d_1 = d(0) , ... , d_{n+1} = d (d_n) , ... $}
\sm
\no We also put
\sm
\ce{$\ell_0 = 0 , \ell_n
= d_n - d_{n-1}$ for $n>0$ and $c_n = \sum_{i \leq n} 2^{i-1} \cdot
\ell_i$ for $n\geq 0$.}
\sm
\no Next, in $W[r]$, define a system $\Sigma = \la \sigma_s ; 
\; s \in \twolom \ra \sub \twolom$ satisfying
$| \sigma_s | = d_{|s| + 1}$ and $s \sub t \Loriar \sigma_s \sub
\sigma_t$ as follows:
\sm
$$\eqalign{\sigma_{\la\ra} =& r \re d_1 \cr
\sigma_{\la 0 \ra} =& \sigma_{\la\ra} \cup r \re [d_1 , d_2 ) \cr
\sigma_{\la 1 \ra} =& \sigma_{\la\ra} \cup r \re [d_2 , d_2 + \ell_2 ) \cr
... &\cr
\sigma_{s_k \ha\la 0 \ra} =& \sigma_{s_k} \cup r \re [ c_{n+1} + 
2 \cdot k \cdot \ell_{n+2} , c_{n+1} + (2k + 1) \cdot \ell_{n+2} ) \cr
\sigma_{s_k \ha \la 1 \ra} =& \sigma_{s_k} \cup r \re [ c_{n+1}
+(2k + 1) \cdot \ell_{n+2} , c_{n+1} + (2k + 2) \cdot \ell_{n+2} ) \cr
}$$
where $\la s_k ; \; k < 2^n \ra$ is the lexicographic enumeration
of sequences of length $n$. Let $T$ be the closure of $\Sigma$
under initial segments; i.e. $T = \{ \sigma_ s \re n ; \; s \in
\twolom \;\land\; n \in\omega \}$. An easy calculation shows
that $T$ must be a perfect tree (otherwise construct a null set
in $W$ which contains $r$, a contradiction).
\par
We proceed to show that $[T]$ is a perfect set of mutually random
reals over $V$. To this end, fix $n\in\omega$, and take a null
set $S \sub (\twoom)^n$ in $V$. We have to prove that
$\{ x = (x_0 , ... , x_{n-1})  \in [T]^n ; \;$ all $x_i$
distinct$\} \cap S = \em$. We start with making some manipulations
with $S$ (most of these arguments are due to T. Bartoszy\'nski,
see [BJ 1] or [Ba]).
\sm

{\capit 2.2. Claim.} {\it We can find a sequence $\la J_m ; \; m\in\omega\ra$
with $J_m \sub (2^m)^n$ such that $\sum_m {|J_m| \over (2^m)^n} < \infty$
and $S \sub \{ x \in (\twoom)^n ; \; \exists^\infty m \; (x \re m \in
J_m ) \}$.} $\qed$
\sm

\no (This is a standard argument.) Define in $V$ a function $f_S
\in\omom$ by putting
$$\eqalign{f_S (0) = &\; 0 \cr
f_S (m+1) = &\; {\rm the \; first} \; k \; {\rm such \; that} \;
(2^{f_S (m)} )^n \cdot \sum_{i=k}^\infty { |J_i| \over (2^i)^n}
< \epsilon_m \cr}$$
where $\la \epsilon_m \in \RR^+ ; \; m\in\omega\ra \in V$ is
strictly decreasing and $\la \epsilon_m \cdot
2^{(2m+1)\cdot n} ; \; m\in\omega\ra$ is summable.
We now get:
\sm

{\capit 2.3. Claim.} {\it $\forall^\infty m \in \omega : \;\;
(2^{d_m})^n \cdot \sum_{i=d_{m+1}}^\infty {|J_i| \over (2^i)^n }
< \epsilon_m$.}
\sm
{\it Proof.} As $d$ eventually dominates $f_S$, we have 
$\forall^\infty m$:
$$(2^{d_m})^n \cdot \sum_{i=d_{m+1}}^\infty {|J_i| \over (2^i)^n}
\leq (2^{d_m})^n \cdot \sum_{i=f_S(d_m)}^\infty {|J_i| \over (2^i)^n}
\leq (2^{f_S (d_m - 1)})^n \cdot \sum_{i=f_S(d_m)}^\infty
{|J_i| \over (2^i)^n} < \epsilon_{d_m - 1} \leq \epsilon_m$$
(because without loss $d_m \leq f_S (d_m -1)$ and $d_m - 1 \geq m$).
$\qed$
\sm

\no In $W$, we construct two new null sets $S^0$ and $S^1$. To this
end define
$$\eqalign{J_k^0 =& \{ s \in (2^{[d_{2k} , d_{2k+2} )})^n ;
\; \exists i \in [d_{2k+1} , d_{2k+2}) \;\exists t \in J_i \;
(t \re [d_{2k} , i ) = s \re [d_{2k} , i )) \} \cr
J_k^1 =& \{ s \in (2^{[d_{2k+1} , d_{2k+3})})^n ; \; \exists
i \in [d_{2k+2} , d_{2k+3}) \; \exists t \in J_i \; (t \re 
[d_{2k+1} , i) = s \re [d_{2k+1} , i)) \} \cr}$$
and then put $S^j = \{ x \in (\twoom)^n ; \; \exists^\infty k
\; (x \re [d_{2k+j} , d_{2k + 2+j}) \in J^j_k ) \}$ for $j \in 2$.
\sm

{\capit 2.4. Claim.} {\it (a) $S \sub S^0 \cup S^1$; \par
(b) $\mu (S^0) = \mu (S^1) =0$. }
\sm
{\it Proof.} (a) is obvious; for (b) note that for almost all $k$ we have (by 2.3)
$${|J^0_k| \over (2^{d_{2k+2} - d_{2k}})^n} \leq
\sum_{i = d_{2k+1}}^{d_{2k+2}-1} (2^{d_{2k}})^n \cdot
{|J_i| \over (2^i)^n} < \epsilon_{2k},$$
and hence $\sum_{k=0}^\infty {|J^0_k| \over (2^{d_{2k+2} - d_{2k}})^n}
< \infty$. Similarly for the $J_k^1$. $\qed$
\sm

\no Using $S^j$ ($j \in 2$) we construct (still in $W$) a null set
$T^j \sub\twoom$ as follows. First let 
$$\eqalign{K_k^j =& \{ s \in 2^{[c_{2k+j} , c_{2k+2+j} )} ; \; \exists t = (t_0 , ... ,
t_{n-1} ) \in J_k^j \;\;  \exists \; {\rm distinct} \; i_0 , ... , i_{n-1} 
\in 2^{2k+j} \cr & \exists \; {\rm distinct} \; j_0 , ... , j_{n-1} \in 
2^{2k+1+j} \;\;  \forall \ell \in n \; \cr & (s \re [c_{2k+j} + i_\ell \cdot
\ell_{2k+1+j} , c_{2k+j} + (i_\ell + 1) \cdot \ell_{2k+1+j} ) =
t_\ell \re [d_{2k+j} , d_{2k+1+j} ) \;\land \cr & s \re [c_{2k+1+j}
+j_\ell \cdot \ell_{2k+2+j} , c_{2k+1+j} + (j_\ell + 1) \cdot
\ell_{2k+2+j}) = t_\ell \re [d_{2k+1+j} , d_{2k+2+j} )) \} \cr}$$
Then put $T^j = \{ x \in \twoom; \; \exists^\infty k \; 
(x \re [c_{2k+j} , c_{2k+2+j} ) \in K_k^j ) \}$ for $j\in 2$.
\sm

{\capit 2.5. Claim.} {\it $\mu (T^0) = \mu (T^1) = 0$.}
\sm
{\it Proof.} Note that for almost all $k$ we have
$${|K_k^0| \over 2^{c_{2k+2} - c_{2k}} } \leq { { (2^{2k} \cdot 
2^{2k+1} )^n \cdot |J_k^0| \cdot 2^{c_{2k+2} - c_{2k} - n \cdot
(\ell_{2k+1} + \ell_{2k+2})} } \over 2^{c_{2k+2} - c_{2k}} }
= { ( 2^{4k+1} )^n \cdot |J^0_k| \over
(2^{d_{2k+2} - d_{2k}} )^n } < \epsilon_{2k} \cdot 2^{(4k+1)\cdot n}$$
(see the proof of 2.4 for this calculation); hence
$\sum_{k=0}^\infty {|K_k^0| \over 2^{c_{2k+2} - c_{2k}} }
< \infty$. Similarly for the $K_k^1$. $\qed$
\sm

\no Thus, $r$ being random over $W$, we know that $r \notin
T^0 \cup T^1$; in particular there is $m\in\omega$ such that
$\forall k \geq m$ we have $r \re [c_{2k} , c_{2k+2} ) \notin
K_k^0$. Now assume that $x = (x_0 , ... , x_{n-1}) \in [T]^n$,
with all branches being distinct; then there is some $m' \in \omega$
such that $\sigma_{s_i} \sub x_i$, where all the $s_i$ have length $m'$
and are pairwise distinct. Let $\hat m = \max \{ m' , m \}$. Then we see
that $\forall k \geq \hat m$:
\sm
\ce{$x \re [d_{2k} , d_{2k+2}) \notin J_k^0$}
\sm
\no (this follows from the definition of $K_k^0$ from $J_k^0$).
Thus $x \notin S^0$. Similarly we see $x \notin S^1$;
hence $x \notin S$, and the proof of Theorem A is complete. $\qed$
\bigskip

{\it 2.6.} We briefly sketch another way of adjoining
a perfect set of mutually random reals. Let us start with fixing some
notation.
$\BB_n$ denotes the random algebra on $(2^\omega)^n$. For $B \in
\BB_n$ (this means that we can think of $B$ as a positive Borel set in
$(\twoom)^n$) and $i<n$, we let $B_i =\{ x ;\; \exists x_0 ,..., 
x_{i-1}, x_{i+1},
..., x_{n-1} \;$ $((x_0,...,x_{i-1},x,x_{i+1},...,x_{n-1})\in B )
\}$, the projection of $B$ on the $i$--th coordinate. Similarly
for $\Gamma \subseteq n$ with $\vert \Gamma \vert = m$, we put $B_\Gamma =
\{ (x_{\Gamma (0)} , ... , x_{\Gamma(m-1)}) ; \; \forall i \in
n\setminus \Gamma \;\exists x_i \; ((x_0, ... , x_{n-1})\in B)
\}$, where $\Gamma (i)$ denotes the $i$--th element of $\Gamma$
for $i < m$.
\par
Elements of $\PP$ are of the form $p =(k,n,T, \{ t_i
;\; i< n \}, B)$ such that $k , n \in \omega$, $T$ is a finite subtree
of $\twolom$
of height $k+1$ with $n$ final nodes all of which are of length
$k$, enumerated as $\{ t_i ; \; i<n\}$, and $B \in \BB_n$ 
satisfies $B_i \subseteq [t_i]$ for $i<n$ (or, equivalently,
$B \subseteq \prod_{i<n} [t_i]$). We put
\smallskip
\centerline{$p=(k^p, n^p, T^p , \{ t^p_i ; \; i<n^p \} , B^p )
\leq (k^q , n^q, T^q, \{ t^q_i ; \; i<n^q \} , B^q )=q$}
\smallskip
\noindent iff $k^p \geq k^q$, $n^p \geq n^q$, $T^p \supseteq T^q$,
and there is a map $j = j^{pq} : n^p \to n^q$ with
$j \re n^q = id \re n^q$ satisfying
\par
\item{(i)} $t_i^p \supseteq t_{j(i)}^q$ for $i\in n$ and \par
\item{(ii)} given any $\Gamma \in
[n^p]^{n^q}$ such that $j \restrict \Gamma$ is one-to-one
(and onto), we have $\tilde j [B^p_\Gamma ]
\subseteq B^q$, where $\tilde j (x_{\Gamma (0)} , ... ,
x_{\Gamma (n^q - 1)} ) = (x_{j^{-1} (0)} ,..., x_{j^{-1} (n^q -1)}
)$ (this means that $\tilde j$ is a homeomorphism of $(\twoom)^{n^q}$ induced
by a permutation of the indices). \par
\no We leave it to the reader to verify that $(\PP , \leq)$
is a $\sigma$--linked p.o. which adds in a canonical way a perfect
set of mutually random reals. It can also be shown that
$\PP$ adjoins a Cohen real (and thus an unbounded real). $\qed$
\bigskip

We conclude this section with a series of questions.
\sm

{\dunh 2.7. Question.} {\it Does $\PP$ add a dominating real?
(Or: does the existence of a perfect set of mutually random
reals over $V\models ZFC$ imply the existence of a dominating real over $V$?)}
\sm

\no We conjecture the answer to be negative. The answer to the
question in parentheses is negative if the adjective ``mutually"
is dropped [BrJ, Theorem 1]. Although $\PP$ does not seem to fit
into the framework of [BrJ, section 1], one may succeed by modifying the
techniques developed there. Closely related is:  
\sm

{\dunh 2.8. Question.} {\it Let $V \subset W$ be models of
$ZFC$. Assume there is a perfect set of random reals over $V$
in $W$. Does this imply the existence of a perfect set of mutually
random reals over $V$ in $W$?}
\sm

\no Note that if {\it random} is replaced by {\it Cohen} then
a positive answer follows immediately from 1.1.
Figuring out whether the p.o. of [BrJ, section 1] adds a perfect
set of mutually random reals may shed new light on both 2.7 and 2.8.
The following is just a reformulation of [BrJ, Question 3]:
\sm

{\dunh 2.9. Question.} {\it Does the existence of a perfect set of
mutually random reals over $V \models ZFC$ imply the existence
of an unbounded real over $V$?}
\sm

\no We say a p.o. $\QQ$ adds a {\it perfect set of mutually
non--constructible reals} if in the extension $V[G]$,
where $G$ is $\QQ$--generic over $V\models ZFC$, there is
a perfect tree $T\sub\twolom$ such that given $n\in\omega$
and distinct $x_0 , ... , x_n \in [T]$, we have $x_n \notin
V[x_0 , ... , x_{n-1}]$. By 1.1 Cohen forcing adjoins
a perfect set of non--constructible reals.
\sm

{\dunh 2.10. Question.} {\it Does adding a random real add
a perfect set of mutually non--constructible reals?}
\Bigskip

{\dunhg $\S$ 3. The Cohen real model}
\Smallskip

Before starting with the proof of Theorem B, we make a few general comments
concerning ``large" free subsets for functions $f: \RR^n \to \RR$.
First, the restriction to Borel or definable functions is
reasonable for it is easy 
to construct, under $CH$, a function
$f : \RR^2 \to \RR$ all free subsets of which are at most
countable. Next notice 
that in case $m = 1$, we can get for
Borel functions a free subset which has positive measure --- and
also one which is non--meager with the Baire property.
(To see this, simply look at the sets $A_\alpha$,
$\alpha \leq \omega$, where $A_\omega = \{ x\in\RR ; \; f(x) = x \}$
and $A_n = \{ x\in\RR ; \; |x-f(x)| > {1\over n} \}$,
all of which are Borel, and note that at least one of these sets 
must be positive.) 
However, this fails in case $m \geq 2$ already
for continuous functions (e.g., for $f(x,y) = |x-y| + x$).
This should motivate somewhat why we are looking for perfect
free subsets.
\par

We shall
show that the statement claimed to be consistent in Theorem B holds in the model
gotten by adding $\omega_1$ Cohen reals to $L$. We denote by
$\CC_{\omega_1}$ the algebra for adding $\omega_1$ Cohen
reals. We start with two well--known lemmata.
\bigskip

{\capit 3.1. Lemma.} {\it Let $\AA$ be a countably generated complete
subalgebra of $\CC_{\omega_1}$. Let $h : \AA \to \CC_{\omega_1}$
be an embedding. Then $h$ can be extended to an automorphism
of $\CC_{\omega_1}$.}
\sm
{\it Proof.} Find $\alpha < \omega_1$ such that $\AA \leq_c \CC_\alpha$ and
$h [\AA] \leq_c \CC_\alpha$. The forcing $\PP$ from $\AA$ to
$\CC_\alpha$ ($\PP '$ from $h[\AA]$ to $\CC_\alpha$) has a countable
dense subset and so is isomorphic to Cohen forcing. Thus $\CC_\alpha$
can be decomposed as $\CC_\alpha = \AA \times \CC = h[\AA] \times
\CC$. Hence we can easily extend $h$ to an automorphism of $\CC_\alpha$,
and then of $\CC_{\omega_1}$. $\qed$
\sm

{\capit 3.2. Lemma.} (Homogeneity) {\it Let $\AA$ be a countably generated complete
subalgebra of $\CC_{\omega_1}$. For any formula $\phi (x)$,
if $\dot x$ is an $\AA$--name, then $\Bool \phi (\dot x ) \Boor_{\CC_{
\omega_1}} \in \AA$.}
\sm
{\it Proof.} This follows from 3.1 in a similar fashion as [Je 1,
Lemma 25.14] follows from [Je 1, Theorem 63].
\par
Let $a = \Bool \phi (\dot x) \Boor_{\CC_{\omega_1}}$. Assume $a \notin
\AA$. Then there is an automorphism $\pi$ of the algebra $\AA '$
generated by $\AA$ and $a$ which fixes $\AA$ and moves $a$
(see [Je 1, Lemma 25.13]). By 3.1 $\pi$ can be extended to an 
automorphism of $\CC_{\omega_1}$. This is a contradiction. $\qed$
\bigskip

{\it 3.3. Proof of Theorem B.} Start with $V=L$.
Let $g : (\omom)^m \to \omom$ be
a $\Delta^1_n$--function in the generic extension
$V[G]$. For $\sigma \in\omlom$ let $\phi_\sigma$ be
a $\Delta_n^1$--formula so that
\sm
\ce{$\phi_\sigma (x_0 , ... , x_{n-1} ) \Loleriar (x_0 , ... , x_{m-1})
\in g^{-1} ([\sigma]).$}
\sm
\no Stepping into some intermediate extension, if necessary, we
may assume that the codes of all $\phi_\sigma$ are in the ground model
$V$. Let $\dot T$ be a Cohen--name for a perfect tree of mutually Cohen
reals (see 1.1). Take $x_0, ... , x_{m-1} , x_m$ branches from
the generic tree $\dot T [G] =T$, and step into the model
$V [x_0 , ... , x_{m-1} ]$. Then:
\sm
\ce{$V [x_0 , ... , x_{m-1}] \models \exists x \; \forall n \;
(\phi_{x \re n} (x_0 , ... , x_{m-1})). \;\;\;\;\; (\star)$}
\sm
\no This is so because
\sm
\ce{$V[G] \models \exists x \; \forall n \; (\phi_{x \re n}
(x_0, ... , x_{m-1}))$}
\sm
\no (as $g$ is a total function in $V[G]$), and all the parameters
of the formula are in $V[x_0 , ... , x_{m-1}]$, and we have
homogeneity (see 3.2). Now choose $x$ as in $(\star)$;
by homogeneity again, we get
\sm
\ce{$V[G] \models g(x_0 , ... , x_{m-1}) = x$.}
\sm
\no In particular, as $x_m \notin V[x_0, ... , x_{m-1}]$,
$V[G] \models g(x_0, ... , x_{m-1}) \neq x_m$. This shows
$\Delta^1_n ({\cal PFS})$.
\par
It is well--known that adding $\omega_1$ many Cohen reals over
$L$ produces a model for $\Delta^1_2 ({\cal B}) + \neg \Sigma^1_2
({\cal B}) + \neg \Delta^1_2 ({\cal L})$ (see, e.g., [JS, section 3] or 
[BJ 2]).
$\qed$
\sm

{\dunh 3.4. Question.} {\it Is $\neg \Delta^1_2 ({\cal B}) + \forall
n \; \Delta^1_n ({\cal PFS})$ consistent with $ZFC$? ---
Is $\neg \Delta^1_2 ({\cal B}) + \neg \Delta^1_2 ({\cal L}) +
\forall n \; \Delta^1_2 ({\cal PFS})$ consistent with $ZFC$?}
\sm

\no Part one of this question is closely related to question 2.10.
More explicitly, a positive answer to 2.10 together with a homogeneity
argument as above should show that adding $\omega_1$ random
reals to $L$ gives a model for $\neg \Delta^1_2 ({\cal B})
+ \forall n \; \Delta^1_n ({\cal PFS})$.
\Bigskip

{\dunhg $\S$ 4. The Sacks real model}
\Smallskip

The model for the statement of Theorem C is gotten by iterating
$\omega_1$ times Sacks forcing $\SS$ with countable support
over the constructible universe $L$.
\par
As before for Cohen and random forcings, we can talk about 
{\it perfect sets of Sacks reals}; more explicitly, given
models $V \subset W$ of $ZFC$, $T$ is a perfect tree of
Sacks reals in $W$ over $V$ if all $x \in [T]$ in $W$ are
$\SS$--generic over $V$. However, we must be more careful
with this notion for it is {\it not} absolute. To see this
first note
\sm
{\capit 4.1. Lemma.} {\it Sacks forcing adds a perfect set
of Sacks reals.}
\sm
{\it Proof.} Sacks [Sa] proved that if $x$ is $\SS$--generic
over $V$, then every new real $y \in V[x] \sem V$ is $\SS$--generic
over $V$ as well. Thus it suffices to show that there is
a perfect set of new reals in $V[x]$.
\par
This is easy to see. Work in $\twoom$. Define $T\sub \twolom$
by $\sigma \in T$ if $\sigma(0) = x(0)$, $\sigma(1+x(0)) = x(1)$, ... ,
$\sigma(n + \sum_{i<n} x(i)) = x(n)$ (where $|\sigma| \geq n+1+
\sum_{i<n} x(i)$). $T$ is perfect because $| x^{-1} (\{ 1 \})
| = \omega$. Note that from $y \in [T]$ we can reconstruct
$x$. $\qed$
\sm
\no However, if we add a further Sacks real $y$ over the model
$V[x]$ of 4.1, then $y$ defines a new branch $y'$ of $T$ from
which we can reconstruct $x$, $T$ and $y$. Thus $y'$ is not
a minimal degree of constructibility over $V$, and therefore
cannot be Sacks--generic over $V$.
\par
We proceed to characterize the statement $\Delta^1_2 ({\cal PFS}_1)
$. We say a set of reals $A \sub \twoom$ is {\it $\SS$--measurable}
if for every perfect tree $T \sub \twolom$ either $ A \cap [T]
$ contains a perfect subset or $[T] \sem A$ contains a perfect
subset.
\sm
{\bolds 4.2. Theorem.} {\it The following are equivalent: \par
\item{(i)} $\forall x \; (L[x] \neq V)$; \par
\item{(ii)} all $\Delta^1_2$--sets of reals are $\SS$--measurable; 
\par
\item{(iii)} all $\Sigma^1_2$--sets of reals are $\SS$--measurable;
\par
\item{(iv)} $\Delta^1_2 ({\cal PFS}_1)$.}
\par\sm
{\it Proof.} We shall show the implications $(i) \Loriar
(iii) \Loriar (ii) \Loriar (iv) \Loriar (i)$.
\par
$(i) \Loriar (iii)$. This is immediate from the
Mansfield--Solovay perfect set theorem [Je 1, Theorem 98].
\par
$(ii) \Loriar (iv)$. Let $f: \twoom \to \twoom$ be a $\Delta^1_2$--function.
For $n\in\omega$, let $A_n = \{ x\in\twoom ; \; d(f(x) , x) > {1 \over n} \}$
where $d$ is the usual metric on $\twoom$, and let $A_\omega =
\{ x\in\twoom ; \; f(x) = x \}$. $A_\omega$ and the $A_n$ are easily
seen to be $\Delta^1_2$--sets of reals. By $\Delta^1_2 - \SS$--measurability
either $A_\omega$ or $\bigcup_n A_n$ must contain a perfect subset.
In the first case we are done; in the second case note that
already some $A_n$ contains a perfect subset $P$ (otherwise use
a fusion argument to show that $\bigcap_n B_n$ contains a perfect
subset where $B_n = (\bigcup_m A_m ) \sem A_n$; however 
$\bigcap_n B_n = \em$, a contradiction). If $P' \sub P$ is perfect
with small enough diameter, we have $f(x) \neq x$ for all $x\in P'$;
thus, $P'$ is a perfect free subset.
\par
$(iv) \Loriar (i)$. This is the boldface version of Mildenberger's
result [Mi, Theorem 1.3]. $\qed$
\sm
It is immediate from Theorem 4.2 ($(i) \Loleriar (iv)$) that
$\Delta^1_2 ({\cal PFS}_1)$ holds after adding $\omega_1$
Sacks reals.
\bigskip
{\it 4.3. Completion of proof of Theorem C.} It remains to 
show that $\Delta^1_2 ({\cal PFS}_2)$ is false in the model $W$
gotten by adding $\omega_1$ Sacks reals over $L$. For this
we use that the degrees of constructibility of the reals
in $W$ are well--ordered of order type $\omega_1$
(Miller [M]).
\par
Let us consider the following three formulae.
$$\eqalign{\phi_1 (x,y) \Loleriar & \exists \alpha \; (y \in L_\alpha
[x] \;\land\; x \notin L_\alpha [y]) \cr
\phi_2 (x,y) \Loleriar & \exists \alpha \; (x \in L_\alpha [y] \;\land\;
y \notin L_\alpha [x]) \cr
\phi_3 (x,y) \Loleriar & \exists \alpha \; (x \in L_\alpha [y]
\;\land\; y \in L_\alpha [x] \;\land\; \forall \beta < \alpha
\; (x\notin L_\beta[y] \;\land\; y \notin L_\beta [x] )) \cr }$$
Using standard arguments it is easy to see that these formulae
are all $\Sigma^1_2$; furthermore they are mutually exclusive;
finally, because the constructibility degrees are well--ordered,
for each $(x,y)$ there is $i \in \{1,2,3\}$ so the $\phi_i
(x,y)$ holds. This entails that the $\phi_i$, $i\in\{ 1,2,3 \}$,
are $\Delta^1_2$ in the iterated Sacks model.
\par
We next define $f: (\omom)^2 \to \omom$; this is done in a recursive
way by cases. Let $x,y \in\omom$.
\sm
\no{\sanse Case 1.} $\phi_1 (x,y)$.
\par
\item{} Take the $L[x]$--minimal perfect tree $T$ such that:
\par
\itemitem{} $x,y \in [T] \;\land\; \forall z \neq x$ before
$y$ in the well--order of $L[x]$ either $z \notin [T]$ or
$f(x,z) \notin [T]$;
\par
\item{} then put $f(x,y) := $ the $L[x]$--minimal branch of $[T]$ different
from $x$ and $y$.
\par\sm
\no{\sanse Case 2.} $\phi_2 (x,y) \;\lor\; \phi_3 (x,y)$.
\par
\item{} Let $f(x,y) = 0$.
\par\sm
\no To see that $f$ has a $\Sigma^1_2$--graph, let
$\phi (x,y,\bar y, T, \bar T, z,\bar z)$ be the conjunction
of the following formulae:
\sm
\item{(i)} $T$ is a perfect tree $\land\; x,y,z \in [T] \;\land\;
z\neq x,y\;\land\; \phi_1 (x,y)$; 
\par
\item{(ii)} $\bar y = \la y_n ; \; n\in\omega\ra , \bar z =\la z_n ;
\; n\in\omega\ra$ are sequences of reals and $\bar T = \la T_n ; \;
n\in\omega\ra$ is a sequence of perfect trees; 
\par
\item{(iii)} for all $n\in\omega$ we have: $y_n <_{L[x]} y$ and 
$y_n \notin [T] \;\lor \; z_n \notin [T] $ and
$y_n \in [T_n]$ and
$(\phi_1(x, y_n) \;\land\; z_n \in [T_n] \;\land\; z_n \neq x,y_n)
\;\lor\; ((\phi_2 (x,y_n) \lor\phi_3(x,y_n)) \;\land\; z_n = 0)$;
\par
\item{(iv)} there is a countable transitive $ZFC$--model $M$
containing $x,y,\bar y,T,\bar T,z,\bar z$ and satisfying:
\par
\itemitem{(a)} $\forall y' <_{L[x]} y \; \exists n \; (y' = y_n)$
and $\forall z' <_{L[x]} z \; (z' \notin [T] \;\lor\; z' = x \;\lor \;
z'=y)$ and for all perfect trees $T' <_{L[x]} T \; (x\notin [T'] \;\lor\;
y \notin [T'] \;\lor\; \exists n \; (y_n , z_n \in [T']))$; 
\par
\itemitem{(b)} for all $n\in\omega$, if $\phi_1 (x,y_n)$, then:
$\forall z' <_{L[x]} z_n \; (z' \notin [T_n] \;\lor\; z' = x \;\lor \;
z'=y_n)$ and for all perfect trees $T' <_{L[x]} T_n \; (x\notin [T'] \;\lor\;
y_n \notin [T'] \;\lor\; \exists m \; (y_m <_{L[x]} y_n \;\land\; 
y_m , z_m \in [T']))$ and $\forall m \; (y_m <_{L[x]} y_n \;\Loriar\;
(y_m \notin [T_n] \;\lor\; z_m \notin [T_n]))$; 
\par\sm
\no $\phi$ is easily seen to be $\Sigma^1_2$, and we have
$f(x,y) = z$ iff either ($\phi_1(x,y)$ and there are $\bar y, T , \bar T,
\bar z$ with $\phi(x,y,\bar y, T,\bar T,z,\bar z)$) or
($(\phi_2(x,y) \;\lor\; \phi_3 (x,y))$ and $z=0$).
\par
We have to check that $f$ doesn't have perfect free subsets.
Suppose $T$ is a perfect tree. Find $x \in [T]$ such that $T \in L[x]$
and $x \notin L[T]$ (this can be done easily by the fact that
the constructibility degrees are well--ordered). Then there
are $\omega_1$ many $y\in L[x]$ with $y \in [T]$
and $x \notin L[y]$; i.e. in particular $\phi_1 (x,y)$
holds. Hence we are in case 1. Because the countably many
predecessors of $T$ in the well--order of $L[x]$ were considered
at some point, we must have $f(x,y) \in [T]$ and
also $f(x,y) \neq x,y$ for one of those $y$. This 
completes the proof of Theorem C. $\qed$
\sm
{\dunh 4.4. Question.} {\it Is $\Delta^1_2({\cal PFS}_2) + \neg
\Delta^1_2 ({\cal PFS}_3)$ consistent with $ZFC$?}

\Bigskip

{\dunhg $\S$ 5. Superperfect trees of Cohen reals}
\Smallskip
We say $T \sub\omlom$ is a {\it superperfect tree of Cohen reals}
over $V$ if $T$ is superperfect and all branches of $T$ are
Cohen over $V$. The next result answers a question addressed
by O. Spinas [Sp, question 3.5].
\sm
{\bolds 5.1. Theorem.} {\it Assume $W \supseteq V$ contains a superperfect
tree of Cohen reals over $V$. Then $W$ contains a real dominating
the reals of $V$.}
\sm
{\it Proof.} Given a real $f\in\omom$, we define the tree $T_f$ by:
$$\sigma \in T_f \Loleriar \forall i \leq |\sigma| \; (f (\max_{j<i}
\sigma (j)) \geq i ).$$
Then $[T_f]$ is nowhere dense (given $\sigma \in T_f$, $\sigma \ha
(0 \re [ |\sigma| , \ell ) ) \notin T_f$ where $\ell > f
( \max_{j < |\sigma|} \sigma (j))$).
\par
Next let $T$ be an arbitrary superperfect tree. Let $\sigma \in
split (T)$. Define $g_T^\sigma \in\omom$ by
$$g_T^\sigma (n) = \min \{ | \tau | ; \; \tau \in split (T) \;\land\;
\exists m \geq n \; (\sigma \ha \la m \ra \sub \tau ) \}.$$
Let $g_T$ be any function eventually dominating all $g_T^\sigma$,
$\sigma \in split (T)$.
\par
To conclude the argument, it suffices to show that $[T_f] \cap [T]
\neq \em$ whenever $g_T$ does not eventually dominate $f$,
$f(i) > | stem (T) |$ for all $i$ and $f$ is increasing. 
To this end, we shall recursively construct $\sigma_n \in split (T)
\cap T_f$ with $\sigma_n \subset \sigma_{n+1}$. $\sigma_0 : = stem (T)$.
To see that $\sigma_0 \in T_f$ use the assumption $\forall i \; (f(i) >
| stem (T) | )$. Assume $\sigma_n \in split (T) \cap T_f$ is
constructed. As $g_T$ eventually dominates $g_T^{\sigma_n}$ and
$f$ is not bounded by $g_T$, there is $m > 0$ so that
$g_T^{\sigma_n} (m) < f(m)$. Going over to a larger $m$, if necessary,
we may assume $\sigma_n \ha \la m \ra \in T$ (this works by
definition of $g_T^{\sigma_n}$ and because $f$ is increasing).
Choose $\tau \in split (T)$ so that $\tau \supseteq \sigma_n \ha
\la m\ra$ and $| \tau | = g_T^{\sigma_n} (m)$. $\tau \in T_f$ is
now easily read off from the definition of $T_f$. $\qed$
\sm
{\capit 5.2. Corollary.} {\it The following are equivalent:
\par
\item{(i)} for every real $x$, there exists a superperfect tree
of Cohen reals over $L[x]$; \par
\item{(ii)} for every real $x$, the union of the meager sets coded
in $L[x]$ is meager; \par
\item{(iii)} $\Sigma^1_2 ({\cal B})$. \par}
\sm
{\it Proof.} The equivalence of (ii) and (iii) is well--known
(see [JS, section 3] or [BJ 2] for a proof).
\par
(i) $\Loriar$ (ii). Fix a real $x$. There is a Cohen real $c$
over $L[x]$, and --- by 5.1. --- a dominating real $d$
over $L[x,c]$. This entails that the union of the meager
sets coded in $L[x]$ is meager in $L[x,c,d]$
(Truss [Tr], see also [BJ 2]). \par
(ii) $\Loriar$ (i). Fix a real $x$. Take a meager set $A$ containing
all meager sets coded in $L[x]$, and construct recursively a superperfect
tree $T$ so that $[T] \cap A = \em$. $\qed$
\bigskip

Our earlier considerations suggest a strengthening of the concept
of a superperfect tree of Cohen reals to the one of a superperfect
tree of {\it mutually} Cohen reals. However, the existence of
such an object is simply inconsistent as is shown by:
\sm
{\capit 5.3. Proposition.} {\it There is a meager set $M \sub (\omom)^2$
such that for all superperfect trees $T$ there are $f \neq g \in [T]$ 
with $\la f,g \ra \in M$.}
\sm
{\it Proof.} Given $\la\tilde\sigma , \tilde\tau\ra \in (\omlom)^2$
with $|\tilde\sigma| = |\tilde\tau|$, we define the tree
$T_{\la \tilde\sigma , \tilde\tau \ra}$ by:
$$\la\sigma,\tau\ra \in T_{\la\tilde\sigma,\tilde\tau\ra} \Loleriar
|\sigma| = |\tau| \;\land\; \tilde\sigma \sub \sigma \;\land\;
\tilde\tau \sub \tau \;\land\; \max_{i<|\sigma|} \{ \tau (i) ,
\sigma (i) \} \geq | \sigma | - | \tilde\sigma |.$$
Then $[T_{\la\tilde\sigma , \tilde\tau\ra} ]$ is nowhere
dense (given $\la\sigma , \tau\ra \in T_{\la \tilde\sigma , \tilde\tau\ra}
$, $\la \sigma , \tau \ra \ha ( 0 \re [ | \sigma | , \ell ))
\notin T_{\la \tilde\sigma , \tilde\tau \ra}$, where
$\ell > \max_{i< |\sigma|} \{ \tau (i) , \sigma (i) \} +
| \tilde\sigma |$). Let $M = \bigcup_{\la \tilde\sigma , \tilde\tau
\ra} [ T_{\la \tilde\sigma , \tilde\tau \ra} ]$.
\par
Given $T \sub \omlom$ superperfect, let $\sigma_0 =\tau_0 = stem (T)$.
We construct recursively $\sigma_n, \tau_n \in split (T)$ with
$\sigma_n \subset \sigma_{n+1}$, $\tau_n \subset \tau_{n+1}$,
$|\sigma_n | < |\tau_n | < |\sigma_{n+1} |$, and $\sigma_1
\not\sub \tau_1$, such that $\la \sigma_n , \tau_n \re | \sigma_n
| \ra \in T_{\la \sigma_0 , \tau_0 \ra}$ and $\la \sigma_{n+1}
\re | \tau_n | , \tau_n \ra \in T_{\la\sigma_0 , \tau_0 \ra}$.
Given $\sigma_n$ and $\tau_n$, choose $i > | \tau_n |$ such that
$\sigma_n \ha \la i \ra \in T$, and let $\sigma_{n+1} \in
split (T)$ with $\sigma_n \ha \la i \ra \sub \sigma_{n+1}$.
Then clearly $\la \sigma_{n+1} \re | \tau_n | , \tau_n \ra \in
T_{\la \sigma_0 , \tau_0 \ra}$. The construction
of $\tau_{n+1}$ is analogous. $\qed$
\bigskip

Of course, a natural problem to investigate in this context
is the existence of superperfect free subsets for Borel functions
$f : (\omom)^n \to \omom$. We first notice that we can always
find such subsets in case $n = 1$. The point is that --- as remarked
already at the beginning of $\S$ 3 ---, given such a function $f$,
there is a non--meager free subset with the property of Baire.
However, a non--meager set with
the property of Baire is easily seen to contain the branches
of a superperfect tree. On the other hand we have the following:
\sm
{\capit 5.4. Proposition.} {\it There is a Borel function
$f : (\omom)^3 \to \omom$ without a superperfect free subset.}
\sm
{\it Proof.} In [Ve, Theorem 2], Veli{\v c}kovi\'c constructed
a partial continuous function $f: D \to \omom$, where
$D \sub (\omom)^3$ is $G_\delta$, such that $f[[T] \cap D]
= \omom$ for any superperfect tree $T \sub \omlom$. Clearly
any Borel extension of $f$ to $(\omom)^3$ satisfies the
requirements of the Proposition. $\qed$
\sm
{\dunh 5.5. Question.} {\it Let $f: (\omom)^2 \to \omom$ be
Borel. Does $f$ necessarily have a superperfect free subset?}
\sm
\no This question might be related to the following:
\sm
{\dunh 5.6. Question.} (Goldstern) {\it Let $\PP$ and $\QQ$
be forcings which both adjoin an unbounded real. Does
$\PP \times \QQ$ add a Cohen real?}
\sm
\no Shelah (unpublished) observed that the product of three p.o.'s
adding an unbounded real adjoins a Cohen real. This is also
a consequence of Veli{\v c}kovi\'c's result mentioned above,
for the real $f(x_0,x_1,x_2)$ is Cohen whenever $x_0, x_1,
x_2$ are independent unbounded reals and $f$ is
Veli{\v c}kovi\'c's function.

\Bigskip

\centerline{\capitg References}
\Smallskip
\itemitem{[Ba]} {\capit T. Bartoszy\'nski,} {\it On covering of real
line by null sets,} Pacific Journal of Mathematics, vol. 131 (1988),
pp. 1-12.
\smallskip
\itemitem{[BJ 1]} {\capit T. Bartoszy\'nski and H. Judah,} {\it
Jumping with random reals,} Annals of Pure and Applied Logic,
vol. 48 (1990), pp. 197-213.
\smallskip
\itemitem{[BJ 2]} {\capit T. Bartoszy\'nski and H. Judah,}
{\it Measure and category: the asymmetry,} forthcoming book.
\smallskip
%\itemitem{[BJS]} {\capit T. Bartoszy\'nski, H. Judah and S.
%Shelah,} {\it The Cicho\'n diagram,} Journal of Symbolic
%Logic, vol. 58 (1993), pp. 401-423.
%\smallskip
\itemitem{[BrJ]} {\capit J. Brendle and H. Judah,} {\it Perfect
sets of random reals,} Israel Journal of Mathematics, vol. 83 (1993),
pp. 153-176.
\smallskip
\itemitem{[BrJS]} {\capit J. Brendle, H. Judah and S. Shelah,}
{\it Combinatorial properties of Hechler forcing,} Annals
of Pure and Applied Logic, vol. 58 (1992), pp. 185-199.
\smallskip
\itemitem{[Je 1]} {\capit T. Jech,} {\it Set theory,} Academic Press,
San Diego, 1978.
\smallskip
\itemitem{[Je 2]} {\capit T. Jech,} {\it Multiple forcing,}
Cambridge University Press, Cambridge, 1986.
\smallskip
\itemitem{[JS]} {\capit H. Judah and S. Shelah,} {\it 
$\Delta_2^1$--sets of reals,} Annals of Pure and Applied Logic,
vol. 42 (1989), pp. 207-223.
\smallskip
\itemitem{[Mi]} {\capit H. Mildenberger,} {\it Perfect free subsets,}
preprint.
\smallskip
\itemitem{[M]} {\capit A. Miller,} {\it Mapping a set of reals
onto the reals,} Journal of Symbolic Logic, vol. 48 (1983), pp.
575-584.
\smallskip
\itemitem{[Mo]} {\capit Y. Moschovakis,} {\it Descriptive Set Theory,}
North--Holland, Amsterdam, 1980.
\smallskip
\itemitem{[My 1]} {\capit J. Mycielski,} {\it Independent sets
in topological algebras,} Fundamenta Mathematicae, vol. 65 (1964),
pp. 139-147.
\smallskip
\itemitem{[My 2]} {\capit J. Mycielski,} {\it Algebraic independence
and measure,} Fundamenta Mathematicae, vol. 71 (1967), pp. 165-169.
\smallskip
\itemitem{[Ox]} {\capit J. C. Oxtoby,} {\it Measure and category,}
Springer, New York Heidelberg Berlin, 2nd edition, 1980.
\smallskip
\itemitem{[Sa]} {\capit G. Sacks,} {\it Forcing with perfect closed sets,}
in: Axiomatic Set Theory, Proc. Symp. Pure Math., vol. 13 (1971),
pp. 331-355.
\smallskip
\itemitem{[Sp]} {\capit O. Spinas,} {\it Generic trees,} preprint.
\smallskip
\itemitem{[Tr]} {\capit J. Truss,} {\it Sets having calibre 
$\aleph_1$,} Logic Colloquium 76, North-Holland, Amsterdam,
1977, pp. 595-612.
\smallskip
\itemitem{[Ve]} {\capit B. Veli{\v c}kovi\'c,} {\it Constructible
reals and analytic sets,} preprint.
\smallskip

\vfill\eject\end